\documentclass[12pt,a4paper,reqno]{amsart}

\usepackage{amsmath,amsthm,amsfonts,amssymb,euscript,cite, mathrsfs,mathtools}
\usepackage{graphics}
\usepackage{xcolor} %
\usepackage{hyperref}
\hypersetup{colorlinks,linkcolor={blue},citecolor={blue},urlcolor={red}}  
\usepackage{enumitem} %
\usepackage{ulem}
\usepackage[noabbrev,capitalize,nameinlink]{cleveref}

\crefname{equation}{}{}
\crefname{algocf}{Algorithm}{Algorithms}
\crefname{equation}{}{} %
\crefname{algocf}{Algorithm}{Algorithms}

\usepackage{graphicx}
\usepackage{color}
\usepackage{transparent}

\usepackage{fullpage} %

\usepackage{seqsplit}

\definecolor{green}{rgb}{0,0.8,0} %

\definecolor{babypink}{rgb}{0.96,0.76,0.76}

\newcommand{\ang}{{\not\negmedspace\partial }}

\newcommand{\rt}{{\mathbb R^3}}

\newcommand{\pa}{\partial}

\newcommand{\gab}{g^{\alpha\beta}}
\newcommand{\hab}{h^{\alpha\beta}}

\newcommand{\pab}{\partial_\beta}
\newcommand{\paa}{\partial_\alpha}

\newcommand{\pat}{\partial_t}
\newcommand{\pav}{\pa_v}

\newcommand{\pao}{\partial\mkern-10mu /\,}

\newcommand{\ti}{\tilde}

\newcommand{\la}{\langle}
\newcommand{\ra}{\rangle}

\newcommand{\ls}{\lesssim}

\newcommand{\inv}{^{-1}}
\newcommand{\invh}{^{-\f12}}

\newcommand{\f}{\frac}

\newcommand{\iy}{\infty}

\newcommand{\ltst}{{L^2(\mathbb{S}^2)}}
\newcommand{\crt}{{C^{ R}_{ T}}}
\newcommand{\cut}{{C^{ U}_{ T}}}

\newcommand{\ctr}{{C^T_R}}

\newcommand{\inte}{{C^{<3T/4}_{ T}}}

\newcommand{\co}{{D_{tr}^{R}}}
\newcommand{\dtr}{{D_{tr}}}

\newcommand{\lflt}{{L^5L^{10}}}

\newcommand{\lt}{{L^2}}
\newcommand{\p}{\phi}

\let\pminus\pm
\renewcommand{\pm}{\phi_{\le m}}
\newcommand{\pmo}{\phi_{\le m+1}}

\newcommand{\pmn}{\p_{\le m+n}}
\newcommand{\nm}{\la t-r\ra} %
\newcommand{\jr}{\la r\ra} %
\newcommand{\ju}{\la u\ra}

\newcommand{\jt}{\la t\ra}
\newcommand{\jv}{\la v\ra}

\newcommand{\jrho}{\la\rho\ra}

\newcommand{\jx}{\la x\ra}

\newcommand{\lr}[1]{\left( #1 \right)}

\let\arXiv\arxiv
\def\doi#1{ {\href{http://dx.doi.org/#1}
   {{\mdseries\ttfamily DOI}}}}

\newtheorem{theorem}{Theorem}[section]
\newtheorem{corollary}[theorem]{Corollary}
\newtheorem{lemma}[theorem]{Lemma}
\newtheorem{proposition}[theorem]{Proposition}
\theoremstyle{definition}
\crefname{claim}{Claim}{Claims}

\newtheorem{definition}[theorem]{Definition}

\theoremstyle{remark}
\newtheorem{remark}[theorem]{Remark}
\theoremstyle{conjecture}

\numberwithin{equation}{section}

\newcommand{\x}{\alpha}
\newcommand{\xb}{\beta}
\newcommand{\xd}{\delta}
	
\newcommand{\xg}{\gamma}
	
\newcommand{\eps}{\epsilon}

\newcommand{\xk}{\kappa}
\newcommand{\xl}{\lambda}
	
\newcommand{\xo}{\omega}
	\newcommand{\xO}{\Omega}

\newcommand{\xs}{\sigma}
	\newcommand{\xS}{\Sigma}

\newcommand{\N}{{\mathbb N}}

\newcommand{\R}{\mathbb R}

\newcommand{\Z}{\mathbb Z}

\newcommand{\calF}{\mathcal F}

\newcommand{\calN}{\mathcal N}

\newcommand{\calR}{\mathcal R}

\newcommand{\tpsi}{{\widetilde{\psi}}}

\usepackage{microtype} %

\begin{document}

\title{Pointwise decay for the energy-critical nonlinear wave equation}

\author{Shi-Zhuo Looi}
\address{Department of Mathematics, University of Kentucky, Lexington, 
  KY  40506}

\begin{abstract}%
This second article in a two-part series (following [arXiv:2105.02865], listed here as \cite{L}) proves optimal pointwise decay rates for the quintic defocusing wave equation with large initial data on nonstationary spacetimes, and both the quintic defocusing and quintic focusing wave equations with small initial data on nonstationary spacetimes. We prove a weighted local energy decay estimate, and use local energy decay and Strichartz estimates on these variable-coefficient backgrounds. By using an iteration scheme, we obtain the optimal pointwise bounds. In addition, we explain how the iteration scheme reaches analogous pointwise bounds for other integral power nonlinearities that are either higher or lower than the quintic power, given the assumption of global existence for those powers (and in the case of the lower powers, given certain initial decay rates).
\end{abstract}

\maketitle

\section{Introduction}

We study the energy-critical nonlinear wave equation in three spatial dimensions on a variety of spacetimes, which we call nonstationary, that are allowed to depend on both $t$ and $x$. For clarity, we emphasise that they are allowed to depend only on $x$ as well. The goal is to obtain the optimal pointwise decay rate, stated in \cref{thm:main}; this is achieved by an iteration scheme that is outlined in \cref{outline of iteration}. Along the way, we prove an $r$-weighted integrated local energy decay estimate (\cref{the rp est}). The results can be viewed as extensions of the results for the linear problem studied in \cite{L}, albeit under stronger assumptions on the coefficients of the wave operator $P$ than in \cite{L}. See also \cref{postthmrem} for how we reach analogous pointwise decay rates for both focusing and defocusing power nonlinearities that are cubic or higher order; the present article focuses on the quintic order case. 

We consider the operator \begin{equation}\label{P def}
P := \paa\gab(t,x)\pab + g^\xo(t,x) \Delta_\xo + B^\x(t,x)\paa + V(t,x) \quad \text{on }\R^{1+3}
\end{equation} where the coefficients are allowed to depend on $t$ and we use the summation convention. %
Here $\Delta_\xo$ denotes the Laplace operator on the unit sphere, and $\x,\xb$ range across $0, \dots, 3$.  The main assumptions on $P$ are that it is hyperbolic and a small asymptotically flat perturbation of the d'Alembertian $\Box = -\pat^2+\Delta$; see \cref{P.assptns} for the precise assumptions on $P$.  
	\,The precise conditions on the potential $V$, the coefficients $B,g^\xo$ and the Lorentzian metric $g$ are given in the main result, \cref{thm:main}. We also remark on the iteration for obtaining pointwise decay rates for other integer powers, both lower and higher than the quintic, in \cref{postthmrem}.

We study the nonlinear Cauchy problem
\begin{equation}   \label{eq:problem}
\begin{cases}
  P\phi(t,x) =  \mu\phi(t,x)^5 & (t,x) \in (0,\infty) \times \R^3\\
(\p(0,x), \pat \p(0,x)) =(\p_0,\p_1)
\end{cases}, \quad \mu \in \{-1, 0 , 1\};
\end{equation} the convention we adopt is that $\mu=1$ corresponds to the defocusing sign. 

	Our main theorem (\cref{thm:main}) states, informally, that if the coefficients of $P - \Box$ are small and asymptotically flat, then the solution to \eqref{eq:problem}, as well as its vector fields, obey the global pointwise decay rates of $\la t-r\ra^{-1 - \min(c(P), 2)}\la t+r\ra^{-1}$; here $c(P)$ is a constant depending on the coefficients in \cref{P def}. The rate of decay coincides with the one obtained by \cite{Gri} in the case $P=\Box$, but in contrast to \cite{Gri} we obtain it, on these general backgrounds $P$, for the solution everywhere in spacetime for initial data in the wider class \cref{data.assptns}, rather than merely in forward light cones $|x| \le \xl t, \xl < 1$ for the Minkowski background and for the narrower class of compactly supported initial data. We believe this pointwise rate of decay to be sharp.
 	\, An overview of the proof is contained in \cref{ss:outline}, where both an outline of the proof and the main novelties of the proof are explained.
	
	\subsection*{History of the problem}
The theory of global existence, uniqueness and scattering for the semilinear wave equation on Minkowski spacetime, in three spatial dimensions, for
$$\Box \phi = \pminus \p^{p+1}, \qquad \p(0,x) = \phi_0(x), \qquad \pa_t \phi(0,x) = \phi_1(x)$$ 
was studied extensively; for instance, in the articles \cite{Jor,Str,Pec2,ShStr,BSh,Gri}. For small initial data, there is a unique global solution if $p>\sqrt{2}$; see \cite{John,GLS,Ta}.
	Work has also been done for the pointwise decay of solutions; see \cite{Pec1,Str,Yang}. In the case of compactly supported smooth data, decay rates were proved in \cite{Sz} (for small data) and in \cite{Gri,BSz} (for large data). 

We now briefly remark on other spacetimes. Much work has been done for solutions to the initial value problem
$$\Box_g \phi = |g|^{-1/2} \pa_\mu(|g|^{1/2} g^{\mu\nu} \pa_\nu\p) = 0, \qquad \p(0,x) = \phi_0(x), \qquad \pa_t \phi(0,x) = \phi_1(x)$$ 
for various Lorentzian metrics $g$. For the Schwarzschild metric, the solution to the wave equation was conjectured to decay at the rate of $t^{-3}$ on a compact region in \cite{Pri}. This rate of decay was shown to hold for the Schwarzschild spacetime, the subextremal Kerr spacetime with $|a|<M$, and other spacetimes; see \cite{DSS,Tat,MTT,Hin,AAG3}. 
We continue this discussion in \cref{ss:LEDests}, in which we elaborate upon the history and utility of local energy decay estimates, and then conclude with a discussion of pointwise estimates and asymptotic behaviour of solutions on both Minkowski backgrounds and perturbations thereof. 
	
\subsection*{Statement of the main theorem}

We state some notation that we use throughout the paper.
	We write $X\ls Y$ to denote $|X| \leq CY$ for an implicit constant $C$ which may vary by line. Similarly, $X \ll Y$ will denote $|X| \le c Y$ for a sufficiently small constant $c>0$. 
	\, In $\R^{1+3}$, we consider %
\[
\partial := (\partial_t, \partial_1,\pa_2,\pa_3), \qquad \Omega := (x^i \partial_j -
x^j \partial_i)_{i,j}, \qquad S := t \partial_t + \sum_{i=1}^3x^i \partial_{i},
\]
which are, respectively, the generators of translations, rotations and scaling. We denote the angular derivatives by $\ang$.
	We set
$$Z := (\pa,\Omega,S) $$ and we define the function class $$S^Z(f)$$ to be the collection of real-valued functions $g$ such that $|Z^J g(t,x)| \ls_J |f|$
whenever $J$ is a multiindex. We will frequently use $f = \jr^\x$ for some real $\x \leq 0$, where $\jr := (1 + |r|^2)^{1/2}$. We also define $S^Z_\text{radial}(f) := \{ g \in S^Z(f) : g \text{ is spherically symmetric}. \}$
	We denote
\begin{align}  \label{vf defn}
\begin{split}
\p_{J} &:= Z^J\p := \partial^i \Omega^j S^k u, \  \text{ if } J = (i,j,k)\\
\pm &:= (\p_{J})_{|J|\le m}, \quad \p_{=m} := (\p_{J})_{|J| = m}.
\end{split}
\end{align}

\subsubsection*{Assumptions on $P$}\label{P.assptns}
Let $h=g-m$, where $m$ denotes the Minkowski metric. Let $\xs \in (0,\infty)$ be real. We make the following assumptions on the coefficients of $P$: 
\begin{equation}\label{coeff.assu}
\begin{split}
 h^{\alpha\beta}, B^{\alpha} \in S^Z(\jr^{-1-\xs}) \\
\pat B^{\alpha}, V\in S^Z(\jr^{-2-\xs}) \\
g^\xo \in S^Z_\text{radial}(\jr^{-2-\xs})
\end{split}
\end{equation}
	In addition, suppose that for a sufficiently small $\eps>0$ we have, for $A_j := \{ 2^j \le |x| \le 2^{j+1} \}$ and an arbitrary interval $I \subset \R_+$,
\begin{align}\label{coeff'}
\begin{split}
\sum_{j \ge0} \sup_{I \times A_j } \jx^2 |\pa^2 \hab| + \jx |\pa \hab| + |\hab| \leq \eps \\
  \sum_{j \ge0} \sup_{I \times A_j }  \jx^2 |\pa B^\x|
  + \jx |B^\x| \leq \eps \\
  \sum_{j \ge0}  \sup_{I \times A_j }  \jx^2 |V| \leq \eps, \quad \sup_{I \times \R^3} |x|^2 |V| \le \eps.
\end{split}
\end{align}
We use the assumptions \cref{coeff'} in order to apply Strichartz estimates on such variable-coefficient backgrounds.

\begin{theorem}[Main theorem]\label{thm:main}Let 
$$ \kappa := \min(\sigma,2)$$
and let $\phi$ solve \cref{eq:problem} with the assumptions \cref{coeff.assu,coeff'}.
Fix $m \in \N$. We assume that for a fixed $N \gg m$,  
\begin{equation}
\label{data.assptns}
\p_0\in \lt(\rt), \ \|\jr^{1/2 + \kappa}\pa\p_{\le N}(0)\|_{\lt(\rt)}<\iy.
\end{equation}
\begin{enumerate}
\item
Suppose in addition that 
$$\|\p_0\|_{H^{N+1}(\rt)} + \|\phi_1\|_{H^N(\rt)} \ll 1.$$
	Then for any $\mu \in \{ -1,0,1\}$, we have for all times $t>0$
\begin{equation}\label{eq:main bound 1}
\sum_{|J|=0}^m |\p_J(t,x)| \ls \f1{\la t+|x|\ra\la t-|x|\ra^{1+\kappa}}, 
\end{equation}

\item\label{itemtwo}
If $\mu \in \{0,1\}$, we have for all times $t>0$
\begin{equation}\label{eq:main bound 2}
\sum_{|J|=0}^m |\p_J(t,x)| \ls \f1{\la t+|x|\ra\la t-|x|\ra^{1+\kappa}}. 
\end{equation}
where the implicit constant is allowed to depend on $\|\p\|_{L^5(\R_+;L^{10}(\rt))}$. 
\end{enumerate}
\end{theorem}

\begin{remark}[Commentary on the main theorem]
Thus, the main theorem states that for $\sigma$ close to 0, the solution decays at the rate $|\phi| \ls \la t+r\ra\inv \nm^{-(1+\sigma)}$. For large $\xs$, the solution decays at the rate $\la t+r\ra\inv \nm^{-2}$. Moreover, we have these rates for vector fields of $\phi$. The parameter $\sigma$ can be of any size (in particular, arbitrarily small positive) and appears in, for instance, the assumptions for proving the $r$-weighted local energy decay estimate of \cref{the rp est}.
\, %
\end{remark}

\begin{remark}[The linear problem: $\mu = 0$ in \cref{eq:problem}]
When $\mu=0$, \cref{thm:main} states that the solution to 
$$P\p=0, \quad \p(0,x) = \p_{0}(x), \ \pat \p(0,x)= \p_{1}(x)$$
obeys the bound
\begin{equation}\label{linear rate}
\pm \ls \jv\inv \ju^{-1-\xs}, \quad v:= t+r, u :=t-r;
\end{equation}
this result is part of the main theorem in the article \cite{L}.
In \cite{L}, all coefficients were allowed to be large perturbations of the Minkowski metric (thus $\eps\in (0,\iy)$ in that setting), and the assumptions \cref{coeff'} were not needed; \cref{coeff'} were brought in to apply Strichartz estimates for the variable-coefficient backgrounds encoded in $P$, and are unneeded for $\mu=0$ due to the absence of the nonlinearity. For instance, in \cite{L}
no assumptions were made on the second derivatives of $\hab$ in order to obtain the final pointwise bound stated in \cref{thm:main}. %
	
	In addition, for $\mu=0$, the main theorem in \cite{L} shows that the final decay rate \cref{linear rate} still holds as long as a weak local energy decay estimate is assumed to hold. As its name implies, assuming such an estimate is a more general assumption than assuming a local energy decay estimate holds. Weak local energy decay estimates are satisfied in a large variety of situations. 
\end{remark}

\begin{remark}[Differing increments]
The argument shown in this paper straightforwardly yields a proof of a more general version of \cref{thm:main} when the decay increments $\sigma$ differ for the coefficients. (In \cref{coeff.assu}, the increments are all assumed to be equal to $\sigma$.) See also the main theorem in the companion article \cite{L}.
\end{remark}

\begin{remark}[Black hole spacetimes]
All the arguments in this paper can be adapted to the exterior of a ball and hence the proofs in this paper can be applied in the case of black hole spacetimes.
\end{remark}

\subsection{Outline of the paper and strategy of the proof} \label{ss:outline}
Here we overview the proof of \cref{thm:main}, first by presenting an outline of the paper and then stating the main novelties of the proof. %
\begin{itemize}
\item
In \cref{sec:notation} we define notation that is used throughout the article.
\item In \cref{sec:fromLEDtoptw}
we connect pointwise bounds to $L^2$ estimates and norms, thereby connecting local energy decay (see \cref{ss:LEDests}) to pointwise bounds. We show how the derivative of $\phi$ decays at a certain better rate than $\phi$ and its vector fields; this improvement depends on the distance from the light cone $\{r=t\}$ at which one is evaluating the pointwise value of the solution. 
\item
In \cref{sec:preliminaries} we rewrite the equation in a way amenable to our pointwise decay iteration scheme. We state and prove lemmas that are used in the scheme to improve the pointwise decay rates of the solution. Finally, we prove an $r$-weighted local energy decay estimate in \cref{the rp est}. To do so, we use Strichartz estimates on such variable-coefficient backgrounds that satisfy the assumptions in \cref{P.assptns} (more precisely, only \cref{coeff'} is assumed\footnote{The Strichartz estimate is written as Theorem 2 in the published version of \cite{MT}, or Theorem 6 in the current arXiv version of \cite{MT}.} in order to use Strichartz estimates). This result of \cref{the rp est} is then used to improve the decay rate of $\phi$ (and its vector fields)---see \cref{rp gain}. 
\item
In \cref{sec:ext} we prove the final decay rate for $\phi$ and its vector fields in the region exterior to the light cone $\{ r=t\}$, that is, in the region $\{ r \geq t \}$. 
\item
In \cref{sec:int} we prove the final decay rate for $\phi$ and its vector fields in the region inside of the light cone, that is, $\{ r \leq t\}$. 

In \cref{rem:HP.iteration} we explain how the iteration applies to other power nonlinearities, reaching pointwise bounds analogous to the one stated in \cref{thm:main}. 
\end{itemize}

\subsubsection*{Main novelties of the proof}\label{ss:strat}

Compared to the linear problem $P\phi = 0$ (see \cite{L}), the present article additionally employs Strichartz estimates and proves an $r$-weighted local energy decay estimate (see \cref{the rp est}). We show that we are able to use these tools for $\p$ and vector fields of $\phi$; see \cref{VFCor,MTlem1,the rp est}. 

In contrast to the linear problem $P\phi=0$ or to power nonlinearities that are sextic or higher order, the initial global decay rate \cref{u/v decay} (which holds for any function that is finite in the $LE^{1}$ local energy norm) alone is not sufficient for reaching the final decay rates \cref{eq:main bound 1,eq:main bound 2}. The purpose of proving \cref{the rp est} is precisely to obtain a slightly better initial decay rate (for both the solution and its vector fields) than \cref{u/v decay}.

We would like to control the local energy norm and the Strichartz norm (i.e. $\lflt$ norm) of not only the solution $\p$, but also its vector fields $\pm$. Prior to proving the $r$-weighted estimate, we must first (see \cref{MTlem1}) control the nonlinearity and its vector fields on the right-hand side of the Strichartz estimate. We partition the time interval $\R_{+}$ into finitely many sub-intervals $I_k$ so that the nonlinearity norm is small on each $I_k$. This smallness enables us to treat the nonlinearity norm perturbatively (that is, we absorb it to the left-hand side of the Strichartz estimate). However, this perturbative argument comes with the cost of the implicit constant in the Strichartz estimate for $\pm$ now being dependent on the Strichartz norm of the solution $\p$; this explains the appended remark in \cref{itemtwo} (the final bound for the large data problem) in \cref{thm:main}.

In proving the $r$-weighted estimate (\cref{the rp est}), for the small data problem the use of \cref{u/v decay} suffices to control the nonlinearity and finish the proof. This is because the bound \cref{u/v decay} then comes with a small factor, and we can immediately treat the nonlinearity perturbatively by bounding four of the five functions in the nonlinearity using \cref{u/v decay}.\footnote{The reader is encouraged to compare this with the analysis in \cref{postthmrem} on the sextic and higher powers and to see why this approach fails for the lower power nonlinearities (cubic and quartic).} 
	In contrast, for the large data problem this factor can be large; nonetheless the goal will still be to treat the nonlinearity perturbatively. To achieve this, we make use of an inductive argument that takes advantage of the defocusing nature of the nonlinearity.   More precisely, even though the defocusing structure is lost upon application of one or more vector fields to the equation, we are able to make use of the zeroth-order $r$-weighted estimate (wherein \textit{no} vector fields have been applied to \cref{eq:problem}) to prove higher order $r$-weighted estimates. Compared to the higher order case, the zeroth-order case controls an additional type of term, namely the nonlinearity term: see \cref{utilised soon}. The estimate \cref{utilised soon} is then used to prove higher order estimates. The idea is that control of some lower order norms allows one to treat the higher order norm perturbatively.%

Once these tools are in hand, we commence the iteration and prove the final decay rate in relatively short order: see \cref{sec:ext,sec:int}. The iteration scheme used in the present article is not dissimilar to that used in \cite{L}; this scheme is outlined in \cref{outline of iteration} and the complete details of the scheme are located in \cref{sec:ext,sec:int}.

\subsection{Local energy decay (LED) estimates}\label{ss:LEDests}
Before stating what an LED estimate is, we define the LED norms. We consider a partition
of $ \R^{3}$ into the dyadic sets $A_R= \{R\leq \la r \ra \leq 2R\}$ for $R > 1$ and $A_{R=1} = \{ r \ls 1 \}$ . 
We define
\begin{equation}
\begin{split}
 \| \phi\|_{LE} &= \sup_R  \| \la r\ra^{-\frac12} \phi\|_{L^2 ([0,\infty) \times A_R)}  \\
 \| \phi\|_{LE[t_0, t_1]} &= \sup_R  \| \la r\ra^{-\frac12} \phi\|_{L^2 ([t_0, t_1] \times A_R)},
\end{split} 
\label{ledef}\end{equation}
its $H^1$ counterpart
\begin{equation}
\begin{split}
  \| \phi\|_{LE^1} &= \| \pa \phi\|_{LE} + \| \la r\ra^{-1} \phi\|_{LE} \\
 \| \phi\|_{LE^1[t_0, t_1]} &= \| \pa \phi\|_{LE[t_0, t_1]} + \| \la r\ra^{-1} \phi\|_{LE[t_0, t_1]},
\end{split}
\end{equation}
as well as the dual norm
\begin{equation}
\begin{split}
 \| f\|_{LE^*} &= \sum_{R\ge1}  \| \la r\ra^{\frac12} f\|_{L^2 ([0,\infty)\times A_R)} \\
 \| f\|_{LE^*[t_0, t_1]} &= \sum_{R\ge1}  \| \la r\ra^{\frac12} f\|_{L^2 ([t_0, t_1] \times A_R)}.
\end{split} 
\label{lesdef}\end{equation}

We have the following scale-invariant estimate on Minkowski backgrounds:
\begin{equation}\label{localenergyflat}
\|\pa \phi\|_{L^{\iy}_t L^2_x} + \| \phi\|_{LE^1}
 \ls \|\pa \phi(0)\|_{L^2} + \|\Box \phi\|_{LE^*+L^1_t L^2_x}
\end{equation}
and a similar estimate involving the $LE^1[t_0, t_1]$ and $LE^*[t_0, t_1]$ norms. This is called a local energy decay estimate, or integrated local energy decay estimate. 
	Morawetz obtained a local energy decay estimate for the Klein-Gordon equation in \cite{M}. Some other work on local energy decay estimates and their applications can be found in, for instance, \cite{Al,KSS,KPV,MS,MT,SmSo}. For local energy decay estimates for small and time-dependent long range
perturbations of the Minkowski space-time, see for instance \cite{Al,MT2,MS} for time dependent
perturbations, and for example \cite{B,BH,SW} for stationary and
nontrapping perturbations. There is a related family of local energy decay estimates for the Schr\"odinger equation.

Moreover, even for large perturbations of the Minkowski metric, if one assumes the absence of trapping then local energy decay estimates can still hold; see for instance \cite{BH, MST}. In addition, even in the presence of sufficiently weak trapping, then estimates similar to these local energy estimates---albeit with a loss of regularity---have been established; see for instance \cite{BCMP,Chr,NZ,WZ}. However, in the presence of sufficiently strong trapping, local energy estimates fail; see \cite{Ral, Sb}. 

We now remark on how local energy decay relates to two types of asymptotic behaviour, namely pointwise decay rates and scattering. Local energy decay in a compact region on an asymptotically flat region implies pointwise decay rates that are related to how rapidly the metric coefficients decay to the Minkowski metric; see, for example, the works \cite{L,MW,Tat, MTT,OS,Mos,AAG1,AAG2,Mor,Hin2, LiT2}. 
	Local energy decay is also involved in proving scattering (another type of asymptotic behaviour) on variable-coefficient backgrounds. In particular, they imply Strichartz estimates on certain variable-coefficient backgrounds, see \cite{MT}. \cite{LT} used local energy decay to prove scattering for the version of the defocusing problem considered in this paper but with only perturbations to the metric, but the argument extends easily to the version of the problem that includes the lower-order terms and angular terms defined in \cref{P def}. 
\subsection*{Pointwise estimates and asymptotic behaviour}
We begin with works studying the Minkowski background. In \cite{S}, pointwise decay estimates were proven for linear wave equations with a source term using the comparison theorem (positivity of the fundamental solution) in $1+3$ dimensions.  In \cite{BCRS}, numerics were shown for the asymptotic behaviour of small spherically symmetric solutions of nonlinear wave equations with a potential which showed that the dominant tail results from a competition between linear and nonlinear effects.

Still on the Minkowski background, in \cite{Yang}, pointwise decay estimates in $\R^{1+3}$ for various ranges of $p$ in the defocusing nonlinearity $|\p|^{p}\p$ were shown given data in a weighted energy space; more precisely, the solution is shown to decay as rapidly as the linear case for $p+1> (1+\sqrt{17})/2$. The paper \cite{Y1} investigated similar questions for $1+d$ dimensions where $d\ge3$. While prior investigations along these lines of questioning used the time decay of $t\mapsto \int_{\R^{1+d}}|\p(t,x)|^{p+2} \,dx$ for $1<p+1<5$ to study pointwise decay estimates and scattering, \cite{Yang} uses the method introduced in \cite{DR} to obtain pointwise decay estimates, via a weighted spacetime energy estimate for $2<p+1<5$. 

We now comment on various works that study other backgrounds, or general backgrounds that include the Minkowski spacetime as a special case. The work \cite{Toh} considers power type nonlinearities with small initial data on Kerr backgrounds. The work \cite{LooTohNullCond} considers the null condition on nonstationary spacetimes similar to those in the present article; it proves global existence and sharp pointwise decay, assuming a local energy decay estimate holds. On these nonstationary spacetimes, the work \cite{Loo22} proves global existence for wave equations with the null condition $P\p = Q(\pa\p,\pa\p)$ assuming a weak LED estimate and for wave equations with cubic and higher order nonlinearities $P\p =\calN_{\ge\text{cubic}}(\pa^{2}\p,\pa\p,\p)$, given small initial data. Under the assumption of global existence (which holds when, for instance, the initial data is small) this work also proves pointwise decay rates for these families of wave equations and for a family $\calF$ of quasilinear wave equations. This family $\calF$ contains, as a special case, the quasilinear wave equations close to Schwarzschild and Kerr spacetimes whose global existence was proved in \cite{LiT} and \cite{LiT2} respectively. 
The upcoming work \cite{LukOh} obtains sharp pointwise asymptotics, given certain assumptions, for a variety of nonlinearities.

\begin{remark}[High and low power nonlinearities] \label{postthmrem} In this remark we explain how the methods of the present article automatically give either partial or complete (if certain known results are assumed) proofs of pointwise decay rates \cref{high powers final} for various other power nonlinearities. Beyond the present remark, we provide more specifics about this in \cref{rem:HP.iteration}, after the iteration has been presented. 

We shall distinguish between the small and large data cases:
\begin{enumerate}
\item
For small initial data, we consider the 
Cauchy problems
$$P\p = \pminus \phi^{p+1}, \quad p \in \Z_{\geq 2}$$
with smooth and compactly supported initial data that is small in a $H^{n+1}\times H^{n}$ norm.
	If $p+1 > 1+\sqrt{2}$ then for sufficiently small and smooth initial data, there exist smooth global solutions. Then the techniques in the present article prove the decay rate
\begin{equation}
\label{high powers final}
|\p(t,x)| \ls \f1{\la t+|x|\ra\la t-|x|\ra^{1+\min(\xs,p-2)}}.
\end{equation}
	Here we simply mention that 
\begin{enumerate}
\item\label{a}
for $p \geq 5$ the estimate found in \cref{the rp est} is unnecessary; instead, the bound \cref{u/v decay} and \cref{conversion,Minkdcyt} alone suffice to reach \cref{high powers final}. Heuristically, this is because if the nonlinearity contains enough decay, then the initial global decay rate from local energy decay alone (see \cref{u/v decay}) suffices to bootstrap the solution toward the sharp decay rate stated in \cref{high powers final}.
\item The case $p=4$ is the subject of the present article (wherein both large and small data are considered).
\item \label{c}
For $p =2,3$, if one had the initial decay rate\footnote{for some sufficiently small $\xd$---more precisely, for $p=2$, we need $2 - 3\xd >1$, while for $p=3$ we need $2 - 4\xd > 1$. This is because of an integration in the radial variable $\rho$ in the backward light cone $\dtr$ (see \cref{D_tr defn}), when bounding $\int_{\dtr} \rho|\p|^{p+1} dA$. More precisely, we are taking the exponent $\text{ex}>1$ in $1/\jrho^{\text{ex}}$. \cref{p2,p3} follow if \cref{rp gain} holds with $\gamma >1$, $\gamma > 1/2$ respectively, and see also \cref{rem:ref}.}
	\begin{align}
\p|_{r < t/2} \ls \jr^{-1-\xd}, \quad \p|_{r \sim t} \ls \jr^{-1+\xd}\ju^{-2\xd}, \qquad p=2 \label{p2}\\
\p |_{r < t/2} \ls \jr^{-3/4-\xd}, \p|_{r \sim t} \ls \jr^{-3/4+\xd}\ju^{-2\xd},\qquad p=3 \label{p3}
\end{align} 
then the iteration also follows \cref{sec:ext,sec:int} nearly verbatim and one reaches the rate \cref{high powers final} from the method in this paper. See \cite{Toh} for a proof of the small data problem on the Kerr background. In the small data case, additional lemmas become available.
\end{enumerate}

\item For large initial data, we consider the (defocusing) Cauchy problems
$$P\p = |\p|^{p}\phi, \quad p \in 2\Z_{\ge1}$$
where we avoid the odd integer values of $p$ because of the issue regarding the smoothness of the modulus of $\phi$ close to the zero set of $\phi$. (The method presented in this article applies vector fields to the nonlinearity.)

The value $p = 4$ is covered in the present article. For higher $p$ values, we remark that global existence has not been explicitly established in the literature. \, 
	If global existence is assumed for the values $p \ge 5$, then the iteration presented in \cref{sec:ext,sec:int} automatically give \cref{high powers final}, with the same remark in \cref{a} above applying here. 
	For $p = 2$, if global existence is assumed, then the remark in \cref{c} above applies here as well.

\end{enumerate}
\end{remark}

\subsection*{Acknowledgements}
I would like to thank Mihai Tohaneanu for bringing the problem to my attention and for stimulating discussions. 

\section{Notation} \label{sec:notation}
We begin by defining dyadic numbers and dyadic conical subregions.
We work only with dyadic numbers that are at least 1. We denote dyadic numbers by capital letters for that variable; for instance, dyadic numbers that form the ranges for radial (resp. temporal and distance from the cone $\{|x|=t\}$) variables will be denoted by $R$ (resp. $T$ and $U$); thus $$R,T, U\ge 1.$$ 
	We choose dyadic integers for $T$ and a power $a$ for $R,U$---thus $R = a^k$ for $k\ge1$--- different from 2 but not much larger than 2, for instance in the interval $(2,5]$, such that for every $j\in\N$, there exists $j'\in\N$ with 
$a^{j'} = \f38 2^j.$
 
\subsubsection*{Dyadic decomposition of spacetime}
We decompose the region $\{r\le t\}$ based on either distance from the cone $\{r=t\}$ or distance from the origin $\{r=0\}$. We fix a dyadic number $T$. The regions $C_T^R, C_T^U$ and $C^T_R$ are where we shall apply Sobolev embedding. 
\begin{align*}
C_T &:= \begin{cases}
\{ (t,x) \in [0,\iy) \times \rt : T \leq t \leq 2T, \ \ r \leq t\} & T>1 \\
\{ (t,x) \in [0,\iy) \times \rt : 0 < t < 2, \ \ r \leq t\} & T=1
\end{cases} \\
C^R_T &:=\begin{cases}
C_T\cap \{R<r<2R\} & R>1\\
C_T\cap \{0 < r < 2\} & R=1
\end{cases}\\
C^U_T &:=\begin{cases}
\{ (t,x) \in [0,\iy) \times \rt  :  T\le t\le 2T\} \cap \{U<|t-r|<2U\} & U>1\\
\{ (t,x) \in [0,\iy) \times \rt  :  T\le t\le 2T\} \cap \{0< |t-r|<2\} & U=1
\end{cases} %
\end{align*}
	As a subregion inside the forward light cone, we define
\begin{equation*}
C_T^{<3T/4} := \bigcup_{R < 3T/8} C_T^R.
\end{equation*}

Now letting $R > T$, we define
\begin{align*}
C^T_R &:= \{ (t,x) \in [0,\iy) \times \rt  :  r \ge t, T \le t\le 2T, R \le r\le 2R, R \le |r-t| \le 2R\}\end{align*}

\subsubsection{The symbols $n$ and $N$} \label{subsec:N}
Throughout the paper the integer $N$ will denote a fixed and sufficiently large positive number, signifying the highest total number of vector fields that will ever be applied to the solution $\p$ to \eqref{eq:problem} in the paper. 

We use the convention that the value of $n$ may vary by line.%

\subsubsection*{Tildes atop sets} If $\xS$ is a dyadic set as defined above, we shall use $\ti\xS$ to indicate a slight enlargement of $\xS$ on $\xS$'s scale. We only perform a finite number of slight enlargements in this paper to dyadic sets. The symbol $\ti\xS$ may vary by line. 

\subsubsection*{The variables $r,u,v$}
If $x =(x^1,x^2,x^3)\in\R^3$, we write 
\begin{align*}
r &:= |x| =\lr{ \sum_{i=1}^3 (x^i)^2 }^{1/2}, \quad u:= t-r,  \quad v := t+r.
\end{align*}

\subsubsection*{Summation of norms}
Recall the subscript notation \cref{vf defn} for vector fields. Let $\| \cdot\|$ be any norm used in this paper. Given any nonnegative integer $N\ge0$, we write $\|g_{\le N}\|$ to denote $\sum_{|J|\le N} \|g_J\|$. For instance, taking the absolute value as an example of the norm, the notation $|\pm(t,x)|$ means
$$|\pm(t,x)| = \sum_{J : |J| \leq m} |\p_J(t,x)|.$$

\begin{definition} \label{D_tr defn}
Let 
$$\calR_1:= \{ R : R < u/8 \}, \quad\calR_2 := \{ R : u/8 < R < v \}, \qquad u > 0.$$

Let $\R_+ :=[0,\iy)$. 
\begin{itemize}
\item
Let $D_{tr}$ denote the backward light cone with apex $(r,t)$
	\[
	D_{tr} := \{ (\rho,s)  \in \R_+^2: -(t+r) \leq s-\rho\leq t-r, \ |t-r| \leq s+\rho \leq t+r\}.
	\] 
When we work with $D_{tr}$ we shall use $(\rho,s)$ as variables, and $D_{tr}^{ R}$ is short for $D_{tr}^{\rho\sim R}$.
\item
For $R>1$, let 
\begin{equation*}
\co:= D_{tr} \cap \{ (\rho,s) : R < \rho<2R\}
\end{equation*}
and let 
\[
D^{R=1}_{tr} := \dtr\cap\{(\rho,s) : \rho<2 \}. \]
\end{itemize}
\end{definition}

\section{From local energy decay to pointwise bounds}\label{sec:fromLEDtoptw}

In this section we will show that local energy decay bounds imply certain slow decay rates for the solution, its vector fields, and its derivatives---see \cref{inptdcExt,derbound}. 

We start with the following pointwise estimate for the second derivative. We shall use it, for instance, when applying \cref{DyadLclsd} to the functions $ w = \pa\pm$ (that is, when we bound the first-order derivatives pointwise); this will be done in \cref{derbound}. 
\begin{lemma}\label{2ndDeBd'}
Assume $\p$ is sufficiently regular. Then for any point $(t,x)$
 \begin{equation}\label{2ndD}
|\pa^2\p_J(t,x)| \ls
\left(\f1\jr+\f1\ju\right)|\pa\p_{\le|J|+1}| +  \left( 1 + \f{t}\ju \right) \jr^{-2}|\p_{\le|J|+2}|
+ \left( 1 + \f{t}\ju \right) |(P\p)_{\le|J|}|.
 \end{equation}
\end{lemma}
A proof of \cref{2ndDeBd'} can be found in \cite{LooTohNullCond}, which in fact proves \cref{2ndDeBd'} under more relaxed assumptions than those made on $P$ in this article. 

By \cref{u/v decay}, \cref{2ndD} immediately implies for solutions to \eqref{eq:problem}: 
\begin{equation}\label{use}
|\pa^2\p_J| \ls
\left(\f1\jr+\f1\ju\right)|\pa\p_{\le|J|+1}| +  \left( 1 + \f{t}\ju \right) \jr^{-2}|\p_{\le|J|+2}|
\end{equation}

\

The primary estimates that let us pass from local energy decay to pointwise bounds are contained in the following lemma. 
\begin{lemma}
\label{DyadLclsd}
Let $w\in C^4$,  $Z_{ij} :=S^i\xO^j$,  $\mu:= \la \min(r,|t-r|) \ra$, and $\calR \in \{\crt,\cut,\ctr\}$. 
	Then we have  \begin{equation}\label{DyadLclsdBd}
 \| w\|_{L^\iy(\calR)} \ls \sum_{i\le 1,j\le 2} \f1{|\calR|^{1/2}}  \lr{  \|Z_{ij} w\|_{L^2(\calR)} +  \|\mu \pa Z_{ij} w\|_{L^2(\calR)} }.
\end{equation}
where we assume $U \le \f38 T$, $R \le \f38 T$ and $R > T$ in the cases $\cut,\crt,\ctr$ respectively, and $|\calR|$ denotes the measure of $\calR$. 
\end{lemma}

The proof of the $\calR\in\{\cut,\crt\}$ case is contained in \cite{MTT} or \cite{L}, while the proof for $\ctr$ can be found in \cite{L}. We only sketch the proof here. One uses exponential coordinates, which results in $\calR$ being transformed into a region of size $O(1)$ in all directions. The result is then proved using the fundamental theorem of calculus for the $s,\rho$ variables and Sobolev embedding for the angular variables. 

The next proposition yields an %
initial global pointwise decay rate for $\phi_J$ under the assumption that the local energy decay norms are finite. We shall improve this rate of decay in future sections (see \cref{sec:ext,sec:int}) for solutions to \eqref{eq:problem}, culminating ultimately in the final pointwise decay rate stated in the main theorem. 

\begin{proposition} \label{inptdcExt}
Let $T$ be fixed and $\p$ be any sufficiently regular function.
There is a fixed positive integer $k$, %
such that for any multi-index $J$ with $|J|\le N - k$, we have:
\begin{align}\label{u/v decay}
\begin{split}
|\p_{J}| \leq \bar C_{|J|} \|\p_{\leq |J|+k}\|_{LE^1[T, 2T]} \ju^{1/2}\jv\inv. %
\end{split}\end{align} %
\end{proposition}
We only sketch the proof here; full details are provided in \cite{L}. One uses \cref{DyadLclsd}, which proves \cref{u/v decay} except in the wave zone. For the wave zone, an extra Hardy-like inequality \cref{con.Hardy} is used to finish the proof. \cref{Hardy} is proven by multiplying by a cutoff function localised to the wave zone; \cite{L} contains a full proof. 
\begin{lemma}\label{Hardy}
If $f\in C^1$, then
\begin{equation}\label{con.Hardy}
\int_{t/2}^{3t/2} \nm^{-2} f(t,x)^2 dx \ls \int_{t/4}^{7t/4}|\pa_r f(t,x)|^2 dx + \f1{t^2} \left(\int_{t/4}^{t/2} f(t,x)^2 dx + \int_{3t/2}^{7t/4} f(t,x)^2 dx\right).
\end{equation}
\end{lemma}

\subsection{Derivative bounds}
The next proposition shows that the first-order derivative (of solutions to \eqref{eq:problem}) decays pointwise faster by a rate of $\min(\jr,\nm)$. It utilises the initial global decay rate \cref{u/v decay}. The estimates in its proof involve the nonlinearity, but for the quintic nonlinearity as in \cref{eq:problem} it turns out that the global bounds \cref{u/v decay} alone already suffice to make the pointwise decay of the first-order derivative similar to the linear case, which is the content of \cref{derbound}. The reader can find details for the linear problem in \cite{L}. \cref{derbound} will be used in the pointwise decay iteration (see later sections, \cref{sec:ext,sec:int})---more precisely, for iterating upon the \textit{linear} components of the equation, namely those having to do with the operator $P - \Box$. By contrast, the nonlinearity in \cref{eq:problem} does not involve any derivatives, so \cref{derbound} will not be involved in the iteration for the nonlinearity. 
\begin{proposition}\label{derbound}
Let $\p$ solve \eqref{eq:problem}, and assume that 
$$\pmn \ls \jr^{-\x}\jt^{-\xb}\ju^{-\eta}$$
for some sufficiently large $n$. We then have
\begin{equation}\label{claim}
\pa\pm \ls \jr^{-\x}\jt^{-\xb}\ju^{-\eta} \mu\inv, \quad \mu := \la \min(r,| t-r |) \ra.
\end{equation}
\end{proposition}

\begin{proof} Let $\calR \in \{ \cut, \crt, C^T_R\}$. 
Given a function $w$, we have
\begin{equation}\label{1stDeBd}
\|\pa w_{\le m}\|_{\lt(\calR)} \ls \| \f{ w_{\le m+n}}{\mu}  \|_{\lt(\ti\calR)} + \|\jr (Pw)_{\le m} \|_{\lt(\ti \calR)}.
\end{equation} (See \cite{L} or \cite{MTT} for a proof of \cref{1stDeBd}.) By the initial global pointwise estimate \cref{u/v decay}, \cref{1stDeBd} with $w = \phi$ yields
\begin{equation}\label{1stDeBd'}
\|\pa\pm\|_{\lt(\calR)} \ls \| \f\pmn\mu \|_{\lt(\ti\calR)}.
\end{equation}

	Recalling \cref{DyadLclsd}, we have 
\begin{align*}
\|\pa\pm\|_{L^\iy(\calR)} 
&\ls |\calR|\invh \sum_Z \|Z\pa \pm\|_{\lt(\calR)} + \|\mu \pa Z \pa \pm \|_{\lt(\calR)} \\
&\ls |\calR|\invh \lr{ \|\pa\pmn\|_{\lt(\calR)} + \|\mu \pa^2 \pmn\|_{\lt(\calR)}	} \\
&\ls |\calR|\invh \lr{ \|\mu\inv \pmn\|_{\lt(\ti\calR)} + \|\mu \pa^2 \pmn\|_{\lt(\calR)}	} \\
&\ls |\calR|\invh \lr{ \|\mu\inv \pmn\|_{\lt(\ti\calR)} + \| \mu \left( \f1\mu|\pa\pmn| + ( 1 + \f{t}\ju ) \jr^{-2}|\pmn| \right) \|_{\lt(\calR)}	} \\
&\ls |\calR|\invh \lr{ \|\mu\inv \pmn\|_{\lt(\ti\calR)} + \|\mu(1 + \f{t}\ju)\jr^{-2} \pmn \|_{\lt(\calR)} }\\
&\ls |\calR|\invh \|\mu\inv \pmn\|_{\lt(\ti\calR)} 
\end{align*}
which follows by \cref{1stDeBd',use}. The final line follows because $\mu^2(1 + t/\ju) \ls \jr^2$. \ Finally, the claim \cref{claim} follows because $\|\mu\inv \pmn\|_{\lt(\ti\calR)} \ls |\calR|^\f12 \|\mu\inv\pmn\|_{L^\iy(\ti\calR)}$. 
\end{proof}

\section{Preliminaries for the iteration, including the $r^{\gamma}$ estimate}\label{sec:preliminaries}

\begin{remark}[The initial data] \label{rem:id}
Let $w := S(t,0)\p[0]$ denote the solution to the free wave equation with initial data $\phi[0]$ at time 0. Then for any $|J| = O_N(1)$,
\begin{equation}\label{Kir}
w_J(t,x) = \f1{|\pa B(x,t)|} \int_{\pa B(x,t)} (\p_0)_J(y) + \nabla_y (\p_0)_J(y) \cdot (y-x) +t (\p_1)_J(y) \, dS(y).
\end{equation}
By \cref{Kir} and the assumptions $\p_0 \in \lt(\R^3)$, 
$$\|\jr^{1/2 + \xk} \pa\p_{\le N}(0)\|_\lt <\iy, \quad \xk = \min( \sigma, 2 )$$
we have
$$ w_J \ls \jv\inv\ju^{-1-\xk} .$$
\end{remark}

\subsection{Summary of the iteration} \label{outline of iteration}
By \cref{rem:id}, we may assume zero initial data in the following iteration. Second, note that it suffices to prove bounds in $|u| \geq 1$, because the desired final decay rate in $|u| < 1$ already holds by \cref{u/v decay}. Third, we distinguish the nonlinearity and the coefficients of $P - \Box$, and for both of these, we apply the fundamental solution. We iterate these two components in lockstep with one another. 

Due to the domain of dependence properties of the wave equation, we shall first complete the iteration in $\{ u < -1 \}$. For the iteration in $\{ u > 1\}$, the decay rates obtained from the fundamental solution are insufficient in the region $\{ r < t/2\}$. To remedy this, we prove \cref{convrsn}. With the new decay rates obtained from \cref{convrsn}, we are then able to obtain new decay rates for the solution and its vector fields. At every step of the iteration, \cref{conversion} is used to turn the decay gained at previous steps into new decay rates. 

\begin{remark}\label{smllirr.xs}
To simplify the iteration, we shall reduce the value of $\sigma$ if necessary to be equal to some positive irrational number less than the original value of $\xs$. %
We do this to avoid the appearance of logarithms in the iterations for $\phi_1$ and $\phi_2$ (see the decomposition \cref{decomp} below). 
We take $0< \xs \ll 1$. %
	\, In the sections spelling out the details of the iteration, namely \cref{sec:ext,sec:int}, we explain how we reach the final decay rate in \cref{thm:main} (wherein the \textit{original} value of $\xs$ is included in the final decay rate). 
\end{remark}

\subsection{Setting up the problem}\label{settingup}

We rewrite \eqref{eq:problem} as
\[
\Box\p = (\Box - P)\p + F = -\pa_\alpha(h^{\alpha\beta}\pa_\beta\p +B^{\alpha}\p) - g^\xo \Delta_\xo \p - (V-\pa_\alpha B^{\alpha})\p + F, \quad F := \p^5
\]

Using the assumptions \eqref{coeff.assu}, we can write this as
$$\Box\p \in \pa \left(S^Z(r^{-1-\sigma}) \p_{\leq 1}\right) + S^Z(r^{-2-\sigma}) \p_{\leq 2} + F$$
	Pick any multiindex $|J| \leq N_1 - 2$. We have after commuting
\begin{equation}\label{first write}
\Box\p_J \in \pa \left(S^Z(r^{-1-\sigma}) \p_{\leq m+1}\right) + S^Z(r^{-2-\sigma}) \p_{\leq m+2} + F_{\le m}
\end{equation}

Due to the derivative gaining only $\ju\inv$ in the wave zone (see \cref{derbound}), we shall perform an additional decomposition as follows.
First, we note that, for any function $w$, 
\begin{equation}\label{D decomp}
\pa w \in S^Z(r^{-1}) w_{\leq 1} + S^Z(1) \pa_t w, \quad r\geq t/2
\end{equation}
which is clear for $\pa_t$ and $\pa_\xo$, while for $\pa_r$ we write
\[
\pa_r = \frac{S}{r}- \frac{t}{r}\pa_t.
\]

Let $\chi_\text{cone}$ be a cutoff adapted to the region $t/2 \le r \le 3t/2$. We now rewrite \eqref{first write} as
\begin{equation}\label{final write}
\Box\p_J \in S^Z(r^{-2-\sigma}) \p_{\leq m+2} + (1- \chi_\text{cone}) \left(S^Z(r^{-1-\sigma}) \pa\p_{\leq m+1}\right) + \pa_t \left(\chi_\text{cone} S^Z(r^{-1-\sigma}) \p_{\leq m+1}\right) + F_{\le m}
\end{equation}
	We now write $\p_J = \sum_{j=1}^3\p_j$ where 
\begin{equation}\label{decomp}
\begin{split}
\Box \p_1 = G_1, \quad G_1 \in S^Z(r^{-2-\sigma}) \p_{\leq m+2} + (1- \chi_\text{cone}) \left(S^Z(r^{-1-\sigma}) \pa\p_{\leq m+1}\right) \\
\Box \p_2 = \pa_t G_2, \quad G_2\in \chi_\text{cone} S^Z(r^{-1-\sigma}) \p_{\leq m+1} \\
\Box \p_3 = F_{\le m} = G_3
\end{split}
\end{equation}
	Henceforth the convention in \cref{subsec:N} will apply to the symbol $n$.

\subsection{Estimates for the fundamental solution}\label{sec:estsfdmt}

We have the following result, which is similar to previous classical results, see for instance \cite{John}, \cite{As}, \cite{STz}, \cite{Sz2}.
\begin{lemma}\label{conversion}
Let $m\ge0$ be an integer and suppose that $\psi : [0,\iy)\times\R^3\to\R$ solves $$\Box\psi (t,x)= g(t,x), \qquad \psi(0) = 0, \quad \pa_t \psi(0) = 0. $$ 
Define
\begin{equation}\label{hdef}
h(t,r) = \sum_{i=0}^2 \|\Omega^i g (t, r\omega)\|_{L^2(\mathbb{S}^2)}
\end{equation}

Assume that %
$$ h(t,r)  \ls \f{1}{ \jr^\x \la v\ra^{\beta} \la u\ra^\eta }, \quad \x \in (2,3) \cup (3,\iy) , \quad \beta\geq 0, \quad \eta\geq -1/2.$$

Define%
\[
\tilde\eta = \left\{ \begin{array}{cc} \eta -2,& \eta<1   \cr -1, & \eta > 1 
  \end{array} \right. .
\]

We then have in both $\{ u > 1\}$, and $\{u < -1\}$ in the case $\x+\xb+\eta>3$:
\begin{equation}\label{lindcy1}
\psi(t, x)\lesssim \frac{1}{\la r\ra\la u\ra^{\alpha+\beta+\tilde\eta-1}}.
\end{equation}

On the other hand, if $\alpha+\beta+\eta < 3$ and $u < -1$, we have
\begin{equation}\label{lindcy2}
\psi(t, x)\lesssim r^{2-(\alpha+\beta+\eta)}. 
\end{equation}

\end{lemma}

\begin{proof}
A detailed proof of \eqref{lindcy1} can be found in Lemma 5.5 of \cite{L} (see also Lemma 6.1 in \cite{Toh}). The idea is to use Sobolev embedding and the positivity of the fundamental solution of $\Box$ to show that
\[
r\psi \lesssim \int_{D_{tr}} \rho h(s,\rho) ds d\rho,
\]
where $D_{tr}$ is the backwards light cone with vertex $(r, t)$, and use \cref{hdef}.

Let us now prove \eqref{lindcy2}. In this case ${D_{tr}} \subset \{r-t\leq u'\leq r+t, \quad r-t\leq \rho\leq r+t\}$ and we obtain, using that $\la u'\ra \lesssim \rho$ and $\rho > t$ in $D_{tr}$:
\[
r\psi \lesssim \int_{r-t}^{r+t} \int_{r-t}^{r+t} \rho^{1-\alpha-\beta}\la u'\ra^{-\eta} d\rho du' \lesssim \int_{r-t}^{r+t} \la u'\ra^{2-(\alpha+\beta+\eta)} du' \lesssim v^{3-(\alpha+\beta+\eta)}
\]
which finishes the proof.
\end{proof}

For the function $\phi_2$ we will use the following result for an inhomogeneity of the form $\pat g$ supported near the cone. The result is similar to \cref{conversion}, except that we gain an extra factor of $\ju$ in the estimate.

\begin{lemma}\label{Minkdcyt}
Let $\psi$ solve 
\begin{equation}\label{Mink2}
\Box \psi = \pa_t g, \qquad \psi(0) = 0, \quad \pa_t \psi(0) = 0,
\end{equation}
where $g$ is supported in $\{\f12 \leq \f{|x|}t \leq \f32 \}$. Let $h$ be as in \eqref{hdef}, and assume that 
\[
|h| + |Sh| + |\Omega h| + \la t-r\ra |\pa h| \ls \frac{1}{\la r\ra^{\alpha}\la u\ra^{\eta}},  \quad 2 < \alpha <3,  \quad \eta\geq -1/2.
\]

Then in $\{ u > 1\}$, and $\{ u < -1\}$ when $\alpha+\eta > 3$
\begin{equation}\label{lindcy1der}
\psi(t, x)\lesssim \frac{1}{\la r\ra\la u\ra^{\alpha+\tilde\eta}}.
\end{equation}

\end{lemma}
\begin{proof}
Let $\tpsi$ solve $\Box \tpsi = g, \tpsi(0)= 0, \pat \tpsi(0)=0.$
In the support of $g$ we have
\[
(t \partial_i + x_i \partial_t) h \lesssim |Sh| + |\Omega h| + \la t-r\ra |\pa_r h|.
\]	
	By \cref{conversion} (with $\beta=0$) applied to $\nabla \tpsi$,  
$\Omega \tpsi$, $S \tpsi$, and the fact
\[
\la u\ra \partial_t \tpsi \lesssim |\nabla \tpsi| +|S\tpsi|+|\Omega \tpsi|  + \sum_i  | (t \partial_i
+ x_i \partial_t) \tpsi|
\]
the claim follows.
\end{proof}

\subsection{$r^\xg$ decay}

\begin{lemma}[Preliminary LED and Strichartz] \label{MTlem} %
Let $I  \in \{ [T_0,T_1], [T_0, \infty) \}$ be an interval, where $T_1 \ge T_0 \ge0$ are real numbers. If for a sufficiently small $\eps>0$,  \cref{coeff'} holds, and $\psi$ is a function, then
\begin{equation}\label{esta}
\|\psi\|_{(LE^1\cap L^5L^{10})(I\times \R^3)} + \| \pa\psi\|_{L^\iy L^2(I\times \R^3)} 
       	\ls  \|\pa\psi(T_0)\|_{L^2} + \|P\psi\|_{(L^1L^2+LE^*)(I\times \R^3)}.
\end{equation}
\end{lemma} 
\begin{remark} We will assume but not prove \cref{MTlem}. The statement of \cref{MTlem} is obtained by combining Theorem 3 in \cite{MT} and Proposition 8 in \cite{MT2}. The assumptions on the operator $P$ in this paper in \cref{coeff'} satisfy the assumptions in both results. 
\end{remark}

\begin{corollary} [\cref{MTlem} with vector fields] \label{VFCor}
For any $m\ge0$
$$\|\pm\|_{(LE^1\cap L^5L^{10})(I\times \R^3)} + \| \pa\pm\|_{L^\iy L^2(I\times \R^3)} 
       	\ls  \|\pa\pm(T_0)\|_{L^2} + \|(P\p)_{\le m}\|_{(L^1L^2)(I\times \R^3)}.$$
\end{corollary}
\begin{proof}
Since $ P\p_J = (P\p)_J + [P,Z^J]\phi $, 
by \cref{esta} we have on any fixed $I$ with left endpoint $T_0$ the following estimate
\begin{align*}
\|\p_J\|_{LE^1 \cap L^5 L^{10}} + \|\pa \p_J\|_{L^\iy L^2}   
       	&\ls  \|\pa\p_J(T_0)\|_{L^2} + \|(P \p)_J\|_{L^1L^2} +\|[P,Z^J]\phi\|_{LE^*}\\
	&\ls \|\pa\p_J(T_0)\|_{L^2} + \|(P \p)_J\|_{L^1L^2} + \eps \|\p_{\le |J|}\|_{LE^1}.
\end{align*} The second line follows from assumptions on $P$. The claim follows for small $\eps$.
\end{proof}

\begin{proposition}[\cref{VFCor} with no nonlinearity] \label{MTlem1}
Let $\phi$ solve \cref{eq:problem}. Assume the hypotheses on $P$ from \cref{coeff'}. For any interval $I  \in \{ [T_0,T_1], [T_0,\iy) \}$, the following estimate holds:
\begin{equation}\label{bd} 
\|\pm\|_{(LE^1\cap L^5L^{10})(I\times \R^3)} + \| \pa\pm\|_{L^\iy L^2(I\times \R^3)} 
       	\leq C(\|\p\|_{\lflt(I\times \R^3)} , m )  \|\pa\pm(T_0)\|_{L^2(\R^3)}
\end{equation}
Since $\|\phi\|_\lflt$ is bounded by the initial data up to an implicit constant, if the initial data are sufficiently small in the energy norm, then the bound \cref{bd} holds %
without any specific dependence on the size of 
$\|\phi\|_{L^5L^{10}}.$
\end{proposition}

\begin{proof}We first remark that 
$$\|\p\|_{L^5(\R_+; L^{10}(\rt))} < \infty$$
which holds by the main theorem in \cite{LT}. 
	For $m \geq 1$, the bounds
$$\|\pm\|_{L^5(\R_+; L^{10}(\rt))} < \infty$$
can be proven by induction, which we now proceed with. 
Suppose that for some integer $m\geq 0$, $\|\pm\|_{L^5(\R_+;L^{10}(\R^3))} < \iy$; we shall show that \cref{bd} holds.

There are intervals $I_0,\dots, I_n$ ``almost-partitioning'' $[0,\iy)$ with
$t_0:=0\in I_0$, 
    and for $0\le j\le n-1$,
$I_j = [t_j,t_{j+1}],$ while $I_n=[t_n,\iy)$,
such that for all $j$, the following norm obeys the bound
\begin{equation}\label{parti.smallness}
\|\pm\|_{L^5(I_j;L^{10}(\rt)) \cap LE^1(I_j\times\R^3)}	\leq 1/(4C_\text{Stri})^{1/4}, \qquad  I_j := [t_j,t_{j+1}]
\end{equation}

\begin{itemize}
\item
By \cref{MTlem}, we obtain
\begin{align*}
\|\pmo\|_{LE^1(I_j)} &+ \|\pmo\|_{L^5(I_j;L^{10}(\R^3))} 
 + \|\pa\pmo\|_{L^\iy(I_j;L^{10}(\R^3))} \\
  &\leq C_\text{Stri} \left( \|\pa\pmo(t_j)\|_{L^2} + \|(P\phi)_{\le m+1}\|_{L^1(I_j;\lt(\rt))} + \eps \|\pmo\|_{LE^1(I_j)} \right). 
\end{align*} 
	(Here, we bounded the commutator in $LE^*(I_j)$. Here $C_\text{Stri}$ is the constant from the Strichartz estimate from \cref{MTlem}, and $\eps$ from \cref{P.assptns} is chosen sufficiently small.)
	This estimate implies the estimate
\begin{align*}
\|\pmo\|_{LE^1(I_j)} &+ \|\pmo\|_{L^5(I_j;L^{10}(\R^3))} 
 + \|\pa\pmo\|_{L^\iy(I_j;L^{2}(\R^3))} \\
  &\le 2C_\text{Stri} \left( \|\pa\pmo(t_j)\|_{L^2} + \|\pm\|^4_{L^5(I_j;L^{10}(\rt))}\|\pmo\|_{L^5(I_j;L^{10}(\rt))} \right).
\end{align*}
\item
By the hypothesis that $\|\pm\|_{L^5(\R_+;L^{10}(\R^3))}$ is finite,
we apply \cref{parti.smallness}. 
Thus
\begin{align}\label{interval}
\begin{split}
\|\pmo\|_{LE^1[t_j,t_{j+1}]} &+ \|\pmo\|_{L^5( [t_j, t_{j+1} ] ;L^{10}(\R^3))} 
 + \|\pa\pmo\|_{L^\iy( [t_j, t_{j+1}] ;L^{10}(\R^3))} \\
  &\le C \|\pa\pmo(t_j)\|_{L^2}, \quad C = C(\|\phi\|_{L^5(\R_+;L^{10}(\R^3))}, m ) \\
  &\le C \|\pa\pmo\|_{L^\iy( [t_{j-1}, t_j] ; \lt(\rt) ) } \\
  &\le C \cdot C \|\pa\pmo(t_{j-1})\|_\lt \\
  & \le C^{j+1} \|\pa\pmo(0)\|_\lt
\end{split}
\end{align}
Here when $j=n$, the interval is understood to read $[t_n, \iy)$. 

\item
By adding these estimates together, 
we obtain
\begin{align*}
\|\pmo\|_{LE^1(\R_+)} + \|\pmo\|_{L^5(\R_+; L^{10}(\R^3))} + \|\pa\pmo\|_{L^\iy(\R_+; L^{10}(\R^3))} 
&\le C^{j+1} \|\pa\pmo(0)\|_\lt
\end{align*}
which holds because
$$\|\pa\pmo\|_{L^\iy(\R_+;L^{10}(\R^3))} \le \sum_{j=0}^{n} \|\pa\pmo\|_{L^\iy(I_j; L^{10}(\R^3))}.$$
\end{itemize}
\end{proof}

 \begin{remark}[Constants can henceforth depend on the $\lflt$ norm of $\phi$]
 Henceforth, we always allow the implicit constant in estimates to depend on $\|\phi\|_{\lflt([0,\iy)\times\R^3)}$. As a consequence of \cref{MTlem1}, we have
\begin{equation}\label{LE1.finite}
\|\pm\|_{LE^1([0,\iy)\times\R^3)} \ls \|\pa\pm(0)\|_{L^2(\R^3)}.
\end{equation}
\end{remark}

In the following theorem and its subsequent application (\cref{rp gain}) to the pointwise decay problem at hand, we have in mind only an arbitrarily small $\gamma>0$. 
\begin{theorem}[The $r^\xg$ estimate] \label{the rp est} Let $\p$ solve \cref{eq:problem}. 
Let $\gamma < 2\sigma, \gamma < 1$ and let the potential $V$ satisfy $V \in S^Z(\eps/r^2)$.
 Let $T_2 > T_1 \geq 0$.
 
 For any integer $m\geq 0$, we have
\begin{equation}\label{rg.est}
A_{\xg,m}+ E^\xg_{\pm}(T_2)\ls_{\|\p\|_{L^5L^{10}}} E^\xg_{\pm}(T_1) + \|\pa \p_{\le m}\|_{LE(T_1,T_2)}^2 + \|\pa^2\phi_{\le m}\|_{LE(T_1,T_2)}^2
\end{equation}
where the $A,E$ norms are:
$$A_{\xg, m} := \int_{T_1}^{T_2}\int_{\R^3} (\pm)^2 r^{\xg-3} + |\bar\partial\pm|^2 r^{\xg - 1} \, dxdt$$
$$E^\xg_{\pm}(T_1) := \| r^{\xg/2} ( \pao\pm,(\pav+\f1{2r})\pm,\f\pm{r} ) (T_1)\|^2_{L^2(\R^3)}, \  \| r^\x (f_1,\dots,f_n)\| := \sum_{j=1}^n\|r^\x f_j \|.$$
\end{theorem}

\begin{proof} 
Fix $m\ge0$. Let $|J| \leq m$. Fix $0 \leq T_1 < T_2$. Let
$$A_{\xg, J} := \int_{T_1}^{T_2}\int_{\R^3} \p_J^2 r^{\xg-3} + |\bar\partial\phi_J|^2 r^{\xg - 1} dxdt.$$
 
Integrating by parts,
\begin{align}\label{Jmain}
	\begin{split}
	\iint_{[T_1,T_2]\times\R^3} &\Box \p_J(r^\gamma\pav \p_J + r^{\gamma-1} \p_J) \,dx\,dt = 
	\iint_{[T_1,T_2]\times\R^3} - \f{\gamma r^{\gamma-1}}2(\pav\p_J)^2 - \f12(2 - \gamma)r^{\gamma-1}|\pao\p_J|^2 \\
	&-\f{\xg(1-\xg)r^{\xg-3}}2\p_J^2 \, dxdt +\int_{\R^3}- r^\xg\left[ \f12|\pa \p_J|^2 + \pa_r\p_J\pat\p_J + \f1{2}\f{\p_J}{r}\pat\p_J\right]_0^T \,dx
	\end{split}	
	\end{align}
\begin{itemize}
\item
We now manipulate the boundary terms to obtain positive definite terms:
We have
\begin{align*}
\int_{\R^3} -r^\gamma \f1{2}\f{\p_J}{r}\pat{\p_J} \,dx 
	&= \int -r^\gamma \f1{2}\f{\p_J}{r}(\pav-\pa_r){\p_J} \,dx  \\
	&= \int -r^\gamma \f1{2}\f{\phi_J}{r}\pav\p_J \, dx + \int_{S^2}\int_0^\iy  r^\xg \f12 \f{\p_J}{r} \pa_r{\p_J} \, r^2 \, drd\xo  \\
 	&=  \int -r^\gamma \f1{2}\f{\phi_J}{r}\pav\p_J\, dx +\int_{S^2}\int_0^\iy  r^\xg \f12 \f{\p_J}{r} \pa_r{\p_J} \, r^2 \, drd\xo \\
 	&= \int -r^\gamma \f1{2}\f{\p_J}{r}\pav\p_J\, dx +\int_{S^2}\int_0^\iy  r^{\gamma+1} \f14 \pa_r{\p_J}^2 \, drd\xo \\
	&= \int -r^\gamma \f1{2}\f{\p_J}{r}\pav\p_J\, dx +-\int_{S^2}\int_0^\iy \f{\gamma+1}4 r^\gamma{\p_J}^2 \, drd\xo \\
	&= \int -r^\gamma \f1{2}\f{\p_J}{r}\pav\p_J\, dx +-\f{\gamma+1}4 \int_\rt r^{\gamma}\f{{\p_J}^2}{r^2} \,dx
\end{align*}

Thus,
\begin{align}\label{bdry}
\begin{split}
- \int_\rt &r^\gamma\left( \f12|\pao\phi_J|^2  + \f12 (\pav\phi_J)^2 + \f{\gamma+1}4 \f{\phi_J^2}{r^2} + \f12 \f{\phi_J}{r}\pav\phi_J\right)_{T_1}^{T_2} \,dx \\
	&=- \int_\rt r^\gamma \left( \f12|\pao\phi_J|^2 + \left[\f12(\pav\phi_J)^2 + \f18\f{\phi_J^2}{r^2} + \f12 \f{\phi_J}{r}\pav\phi_J \right] + ( \f\gamma4 + \f18 ) \f{\phi_J^2}{r^2}  \right)_{T_1}^{T_2} \,dx \\
	&=-\int_\rt r^\gamma \left( \f12|\pao\phi_J|^2 + \f12\left[ \pav\phi_J + \f{\phi_J}{2r} \right]^2 + ( \f\gamma4 + \f18 ) \f{\phi_J^2}{r^2}  \right)_{T_1}^{T_2} \,dx
\end{split}
\end{align}

\item We shall now prove by induction the claim that
$$\int_{\R_{+}}\int_{\rt}(\pm)^{2}r^{\xg-3}dxdt \le C$$
where $C$ is a constant depending on the initial data.
We make use of the defocusing sign of the nonlinearity. 
We have 
\begin{equation*}
 \iint \phi^5(r^\gamma\pav\phi+r^{\gamma-1}\phi) \,dxdt = \int_\rt r^\gamma\f{\phi^6}6|_0^T \,dx +\iint \left(\f23 - \f\gamma6\right)r^{\gamma-1}\phi^6 \,dxdt
\end{equation*}Note that both terms are nonnegative for our range of small $\xg$ (indeed any $\gamma < 4$).
This implies (for our range of small $\gamma$)
\begin{equation}
\label{utilised soon}\int_0^\iy\int_\rt r^{\xg-1}\p^6 \,dxdt < \iy
\end{equation}
(because $T_{1}$ and $T_{2}$ were arbitrary); and it also establishes the claim for the base case value $m = 0$.

Suppose that for some $m \geq 0$, 
$$\int_0^\iy\int_\rt (\pm)^2 r^{\xg - 3} dxdt < \iy.$$
Then, letting $M \p_J = r^\xg ( \pav \p_J + r^{-1}\p_J ) $ and $|J| = m+1$, 
\begin{align}\label{first}
\int_{T_{1}}^{T_{2}}\int_{\R^3} (\phi^5)_J M\phi_J dxdt 
  &= \int r^{(\xg - 1)/2} M\phi_J \cdot r^{(\xg + 1)/2} (\phi^5)_J dxdt \\
  &\ls \eps' \int r^{\xg - 1} \Big( (\pav\phi_J)^2 + (r\inv \phi_J)^2 \Big) dxdt + \f1{\eps'} \int r^{\xg+1} ( (\phi^5)_J )^2 dxdt
\end{align}for some small $\eps'>0$. We treat that term perturbatively and absorb it to the left hand side. Then we note that:
\begin{itemize}
\item
In the case when there is no single factor in the nonlinearity that has $m+1$ vector fields falling on it, we have by \cref{u/v decay}
$$\int r^{\xg+1} ( (\phi^5)_J )^2 dxdt \ls \int r^{\xg - 3} A^{8} (\pm)^{2} dxdt, \quad A:= \|\pmn\|_{LE^{1}(\R_{+})}$$
which is bounded by a constant depending on the initial data by the induction hypothesis. 
\item
If there is a factor with $m+1$ fields falling on it, we return to \cref{first} and use \cref{utilised soon}. 
$$\int\p^{4}r^{\xg-1} \p_{J}^{2} dxdt \le \eps' \int r^{\xg-2}\p_{J}^{4} + \f1{\eps'} \int r^{\xg}\p^{8} \le \eps' \int r^{\xg-3}A^{2}\p_{J}^{2} + \f1{\eps'} A^{2} \int r^{\xg - 1}\p^{6}$$
for some small $\eps'>0$. 
	We again used \cref{u/v decay}.
	We absorb the small term to the left hand side, and the other term is bounded by a constant depending on the initial data, from \cref{utilised soon}.
\end{itemize}

    \item 
   \begin{enumerate}
    \item Here, in dealing with the potential $V$, we assume only that $V\in S^Z(\eps r^{-2})$. 
    We have: 
    \begin{align}\label{potential}
\begin{split}
\int_{T_1}^{T_2} \int_{\R^3} |V_{\le m} \pm r^\xg (\pav \p_J + \f{\p_J}{r}) | \,dxdt 
   &\ls \eps \int \f1{r^{2-\xg}}|\pm| ( |\pav \p_J|+|\f{\p_J}r| ) \,dxdt \\
   &\ls \eps \int \f{(\pm)^2 + \p_J^2}{r^{3-\xg}}+ \f{|\pav\p_J|^2}{r^{1-\xg}}  \, dxdt \\
   &\ls \eps A_{\xg, m} %
\end{split}
\end{align}
If $B \in S^Z(\jr^{-1-\xs_B})$ and $2 \sigma_B > \gamma$:
\begin{align}\label{B term}
\begin{split}
\int_{T_1}^{T_2} \int_{\R^3} |B_{\le m} \pa\pm r^{\xg-1}\p_J| dxdt
  &\ls  \int \f1{\jr^{1+\xs_B}} |\pa\pm r^{\xg-1}\p_J|  \\
  &\ls \f1\eps\int \f{|\pa\pm|^2}{\jr^{1 + 2\xs_B - \xg}} + \eps \int r^{\xg-3}\p_J^2 \\
  &\ls \f1\eps \|\pm\|_{LE^1(T_1,T_2)}^2 +  \eps A_{\xg, J}
\end{split}
\end{align}    The bound on $\int |B_{\le m}\pa\pm r^\xg| \cdot |\bar\pa\phi_J|dxdt$ is similar.

\item We consider now all of the terms that involve the metric $\hab$. We may schematically write this as $\int (|\pa h_{\le m} \pa\pm| + |h_{\le m} \pa^2\pm|) r^\xg(\f{|\p_J|}r+|\pav\p_J|)dxdt$, where $|J| = m$. 
\begin{align}\label{h term}
\begin{split}
 \int_{T_1}^{T_2} \int_{\R^3} &(|\pa h_{\le m} \pa\pm| + |h_{\le m} \pa^2\pm|) r^\xg(\f{|\p_J|}r+|\pav\p_J|)dxdt \\
	&\ls \eps \int \f1{\jr^{1+\xs}}r^\xg\left(|\pa\pm| + |\pa^2\pm|\right) (|r\inv \p_J| + |\pav\p_J|) \,dxdt\\
	&\ls \eps \int \f{r^{\f{\xg+1}2}}{\jr^{1+\xs}} \left(|\pa\pm| + |\pa^2\pm|\right) \cdot r^{\f{\xg-1}2} (|r\inv \p_J| + |\pav\p_J|) \,dxdt    \\
	&\ls\eps \int \f1{\jr^{2+2\xs}}r^{\xg+1}  \left(|\pa\pm| + |\pa^2\pm|\right)^2 \,dxdt +\eps A_{\xg, J} \\
	&\ls \eps \|\pa\pm\|_{LE(T_1,T_2)}^2 + \eps \|\pa^2\pm\|_{LE(T_1,T_2)}^2 + \eps A_{\xg, J}	    \text{ if } 2 \xs > \xg
\end{split}
\end{align} 
which we can control using \cref{LE1.finite} if we assume $2\xs > \xg$.

    \end{enumerate}
    
\end{itemize}

Taking the sum of \cref{Jmain,bdry,potential,B term,h term} over all $|J| \leq m$, i.e. $\sum_{|J| \leq m}$ (\cref{Jmain,bdry,potential,B term,h term}), 
  and taking into account our argument for the nonlinear terms as well,
  we get
$$A_{\xg, m} + E^\xg_{\pm}(T) \ls E^\xg_{\pm}(0) + \|\pm\|^2_{LE^1(T_1,T_2)} + \|\pa^2 \pm\|_{LE(T_1,T_2)}^2.$$
\end{proof}

Thus we used \cref{P.assptns}'s assumptions to obtain finiteness of the local energy norms $\|\pm \|_{LE^1(\R_+ \times\R^3)}$, and now we have, in particular, showed that
    $$\int_0^\iy \int_{\R^3} (\pm)^2r^{\xg -3} dxdt \ls E^\xg_{\pm}(0) + \|\pa\pmo(0)\|_\lt^2 \le C_0$$
where $C_0$ is a constant depending only on initial data.

    \begin{lemma}\label{aux}
Recall \cref{D_tr defn}. 
Let $v_+ := \la s + \rho \ra$ where $(\rho,s)\in \dtr$. Then 
\begin{equation}\label{R_1 bd}
\|v_+\inv\|_{L^2(D_{tr}^{R\in \calR_1})} \ls 1,
\end{equation}
\begin{equation}\label{R_2 bd}
\|v_+\inv\|_{L^2(D_{tr}^{R\in \calR_2})} \ls \lr{ \f\ju{R} }^\f12.
\end{equation}
\end{lemma}
\begin{proof}
\cref{R_1 bd,R_2 bd} follow from \cref{conversion}.
\end{proof}

\begin{proposition}[Application of the $r^\gamma$ estimate] \label{rp gain}
Let $\p$ solve $P\phi = \p^5$.
Assume the hypotheses on $\gamma$ in \cref{the rp est} and also \cref{coeff.assu}. 
\begin{equation}\label{int.est}
|\jr (\p_3)_{\le m}| \ls \ju^{1/2 -\gamma/2}, \quad u > 1
\end{equation}
\begin{equation}
|(\p_3)_{\le m} | \ls r^{-\f12 -\f\xg2}, \quad u < -1
\end{equation}
\end{proposition}

\begin{remark}\label{rem:ref}
If $\xg \geq 1$, then for $u<-1$ this theorem would instead conclude 
$$(\p_{3})_{\le m}\ls r\inv \ju^{-\f12(\xg-1)}.$$
\end{remark}

\begin{proof} 
\begin{itemize}
\item Let $u > 1$. 
We now show 
$$\int_\dtr \rho H_3 \, dA \ls \ju^{1/2 - \xg/2}, \quad H_3(t,r) := \sum_{k=0}^2 \| \xO^k (\p^5)_{\le m}(t,r\xo)\|_\ltst.$$ 
	
We have\begin{align}\label{hi}
\begin{split}
\int_\co \rho H_3 \, dA
	&\ls \int_\co \rho\f{1}{\la s+\rho\ra^{2}} \| \pmn \|_\ltst \, dA\\
	&\ls\int_\co \f1{v_+} \|\pmn\|_\ltst dA\\
	&\ls \|v_+\inv\|_{\lt(D_{tr}^R)} \|\pmn\|_{L^2_{\rho,s,\xo}}	\\
	&\ls \|v_+\inv\|_{\lt(D_{tr}^R)} \cdot \f1R \f1{R^{(\xg - 3)/2}} C_0
\end{split}
\end{align}where $C_0$ is a constant depending on the initial data.
The first line follows by \cref{u/v decay} and the last line follows by \cref{the rp est}. 

\begin{enumerate}
\item Let RHS denote ``the right-hand side of.''
By \cref{aux}, 
$$\sum_{R\in\calR_1} \text{RHS}\cref{hi} \ls \ju^{1/2 - \xg/2}, \quad \gamma \in (0,1).$$ %
\item
Fix $R\in \calR_2$. By \cref{aux} we have
$$\text{RHS}\cref{hi} \ls \left( \f\ju{R} \right)^\f12 \f1R \f1{R^{(\xg - 3)/2}} 
	= \ju^{1/2} R^{-\xg/2}.$$
Then we have
$$\sum_{R\in\calR_2} \text{RHS}\cref{hi} \ls \ju^{1/2 - \xg/2},  \quad \text{valid for } \gamma \in (0,\infty).$$
\end{enumerate}
This finishes the proof of \cref{int.est}.
\item Let $u < -1$.
\begin{align*}
\int_\dtr \rho H_3 \, dA 
 &\ls \int_\dtr \rho^{\f{\xg - 1}2} \| \pmn \|_{L^2_\xo} \cdot \rho^{\f{3-\xg}2} \| \pmn \|^4_{L^\iy_\xo} \, dA \\
 &\ls \left( \int \rho^{3-\xg}\|\pmn \|_{L^\iy_\xo}^8 dA \right)^\f12\\
 &\ls \left( r^{1-\xg} \right)^\f12.
\end{align*}
The second line follows from \cref{the rp est}. 
The third line follows from \cref{u/v decay}, \cref{conversion} and the assumption $\gamma < 1$. 

\end{itemize}
\end{proof}

\section{The iteration in $\{ u <  -1\}$} \label{sec:ext}

\begin{theorem}\label{ext.thm}
If $u < -1$, then
$$\pm \ls \jr\inv \ju^{-1 - \min(\xs,2)},$$
Here $\xs$ denotes the original value of $\xs$ taken from \cref{thm:main}.
\end{theorem}

\begin{proof}

We begin with the bounds in \cref{u/v decay,derbound,rp gain}, which in the outside region translate to
\begin{equation}\label{inbd}
  \p_{\le m+n} \ls \frac{\ju^{1/2}}{\la r\ra}, \quad \pa\p_{\le m+n} \ls \frac{1}{\la r\ra \la u\ra^{1/2}}, \quad \pmn \ls \f{\ju^{1/2 - \xg/2}}{\jr}.
  \end{equation}
 For simplicity, we shall use the first (far left) $\pmn$ bound for $\p_1$ and the other (far right) $\pmn$ bound for $\p_3$. 
Since $\la u\ra\leq\la r\ra$, \cref{inbd} can be weakened to
\begin{equation}\label{1stbd}
  \p_{\le m+n} \ls \frac{1}{\la r\ra^{1/2}}, \quad \pa\p_{\le m+n} \ls \frac{1}{\la r\ra^{1/2} \la u\ra}, \quad \pmn \ls \frac{1}{\la r\ra^{1/2 + \xg/2}}
\end{equation}

Recall the decomposition \eqref{decomp}, and let 
\[
H_i = \sum_{k=0}^2 \|\Omega^k (G_i)_{\leq n} (t, r\omega)\|_{L^2(\mathbb{S}^2)}.
\]

Let $\xs$ denote the reduced, irrational number mentioned in \cref{smllirr.xs} until stated otherwise. We thus have, using \eqref{1stbd}:
\[
H_1 \ls \frac{1}{\la r\ra^{5/2+\sigma}}, \quad \pat H_2 \ls \frac{1}{\la r\ra^{3/2+\sigma}\la u\ra}, \quad H_3 \ls \f1{\jr^{5/2 + 5/2\xg}}.
\]

By \eqref{lindcy2} with $\alpha=5/2+\sigma$, $\beta = 0$, and $\eta=0$, we obtain
\[
(\p_1)_{\leq m+n} \lesssim r^{-1/2-\xs} %
\]
which gains a factor of $\la r\ra^{-\sigma}$ compared to \eqref{1stbd}. Similarly \eqref{lindcy2} with $\alpha=3/2+\sigma$, $\beta = 0$, and $\eta=1$ yields
\[
(\p_2)_{\leq m+n} \lesssim r^{-1/2-\sigma} %
\]
Finaly, \cref{lindcy2} with $\x = 2 + 2 \xg, \xb=0, \eta = 1/2$ yields
\[
(\p_3)_{\le m+n} \ls r^{-1/2 - 2 \xg}.
\]

The three inequalities above, combined with \cref{derbound}, give the following improved bounds (by a factor of $\la r\ra^{-\xs'}$ where $\xs' := \min( 2\xg , \xs)$).
\begin{equation}\label{2ndbd}
  \p_{\le m+n} \ls \f{1}{\jr^{1/2+\xs'}}, \quad \pa\p_{\le m+n} \ls \f{1}{\jr^{1/2+\xs'} \la u\ra}.
\end{equation}

We now repeat the iteration, replacing $\alpha$ by $\alpha+\sigma'$ and applying \eqref{lindcy2}. The process stops after $\lfloor\frac{1}{2\sigma'}\rfloor$ steps, when \eqref{lindcy2}, combined with \cref{derbound} yield
\begin{equation}\label{3rdbd}
  \pmn \ls \frac{1}{\la r\ra}, \quad \pa\p_{\le m+n} \ls \frac{1}{\la r\ra \la u\ra}.
\end{equation}

We now switch to using \eqref{lindcy1} for $\p_1$ and $\p_3$, and \eqref{lindcy1der} for $\p_2$. Note that \eqref{3rdbd} implies
\[
H_1 \ls \frac{1}{\la r\ra^{3+\sigma}}, \quad H_2 \lesssim \frac{1}{\la r\ra^{2+\sigma}}, \quad H_3\ls \f1{\jr^5}.
\]

By \eqref{lindcy1} with $\alpha=2+\sigma$, $\beta = 1$, and $\eta=0$, we obtain
\[
(\p_1)_{\leq m+n} \lesssim r^{-1} \la u\ra^{-\sigma} 
\]
Similarly \eqref{lindcy1der} with $\alpha=2+\sigma$, and $\eta=0$ yields
\[
(\p_2)_{\leq m+n} \lesssim r^{-1} \la u\ra^{-\sigma} 
\]
Finally, \eqref{lindcy1} with $\alpha=5$, $\beta = 0$, and $\eta=0$ yields
\[
(\p_3)_{\leq m+n} \lesssim r^{-1} \la u\ra^{-2} 
\]

 We now repeat the iteration. We can continue improving the decay rates of $\p_1$ and $\p_2$ all the way to
\begin{equation}\label{final1}
(\p_1)_{\leq m}, (\p_2)_{\leq m} \ls r^{-1} \la u\ra^{-1 - \sigma}.
\end{equation}
For $\p_3$, we note that after the bounds $(\p_1)_{\le m+n}, (\p_2)_{\le m+n} \ls r\inv \ju^{-1/5-}$ are obtained, we have 
\begin{equation}\label{final2}
(\p_3)_{\le m+n} \ls r\inv \ju^{-3}.
\end{equation}
By \cref{final1,final2} we now have, for the \textit{original} value of $\xs$ from \cref{thm:main},
$$ H_1 \ls \f1{\jr^{3 + \xs}\ju^{1+}}, \quad H_2 \ls \f1{\jr^{2 + \xs}\ju^{1+}}, \quad H_3 \ls \f1{\jr^5\ju^{1+}}.$$
Using \cref{lindcy1,lindcy1der} now completes the proof.
\end{proof}

\section{The iteration in $\{ u > 1\}$} \label{sec:int}

\subsection{Converting $r$ decay to $t$ decay}

The pointwise decay rates for the solution and its vector fields obtained from the estimates for the fundamental solution (see \cref{sec:estsfdmt}) are by themselves insufficient for completing the iteration. We show below that if the $\jr\inv$ decay from the fundamental solution is converted into $\jt\inv$, then the iteration does work.  

\begin{lemma}\label{lem:aux}
$$\|\pm\|_{LE^1(\inte)} \ls T\inv \|\jr \pmn\|_{LE^1(\inte)} + \|(\p^5)_{\le m+n}\|_{LE^*(\inte)}.$$
\end{lemma}
We only sketch the proof of \cref{lem:aux} here, and refer the reader to \cite{L} or \cite{MTT} for details. We may assume that $\phi$ is supported in $\inte$.  By the local energy decay estimate, we control $\|\pm\|_{LE^1(\inte)}$ by an energy term plus $\|(\p^5)_{\le m}\|_{LE^*(\inte)}$.  
Then averaging in time is used to bound the energy term by the $LE^1$ norm, leading to the desired estimate.

The next proposition uses \cref{lem:aux} to obtain better pointwise decay for the solution and its vector fields in the region $\{ r < t/2\}$. 

\begin{proposition}\label{convrsn}
Let $\phi$ solve \eqref{eq:problem}. Let $\xd > 0$. Assume that 
\begin{equation}\label{r bds}
\phi_{\le M}|_{r \le 3t/4} \ls \jr^{-1}\ju^{1/2 - n\delta}, \quad \phi_{\le M}|_{r \le 3t/4} \ls \jt\inv \ju^{1/2 - (n-1)\delta}, \quad n \geq 1
\end{equation}
for an $M$ that is sufficiently larger than $m$.
	Then we have 
	$$ \pm|_\inte\ls  \jt\inv \ju^{1/2 - n\delta} .$$
\end{proposition}

\begin{proof}
By \cref{derbound} and \cref{r bds},
$$ T\inv \|\jr \pmn\|_{LE^1} \ls T\inv \|\pmn\|_{LE} \ls T^{ - q}, \quad q = n \delta.$$

 Fix $n \geq 1$.
 For $A_{R = 1}$, we use the latter bound in \cref{r bds} to obtain 
 $$\|(\p^5)_{\le m+n}\|_{L^2[T,2T]\lt(A_{R = 1} )} \ls T^{-2-5(n-1)\xd}.$$
 	 For $A_R, R>1$, we use the former bound and the latter bound in \cref{r bds} in a three-to-two ratio (respectively). This yields 
 $$\|\jr^\f12 (\p^5)_{\le m+n}\|_{L^2[T,2T]\lt(A_R)} \ls \left( T^{-10n\delta + 6 \delta} \right)^\f12 = T^{-5n\delta + 3\delta}, \quad R>1.$$
Therefore \cref{lem:aux} implies, after the dyadic sum, 
$$\|\pm\|_{LE^1(\inte)} \ls T^{- n\delta}$$
and the conclusion now follows by \cref{DyadLclsd}. 
\end{proof}

\subsection{The iteration}

\begin{theorem}
If $u > 1$, then 
$$\pm \ls \jv\inv \ju^{-1 - \min(\xs,2)}.$$
Here $\xs$ denotes the original value of $\xs$ taken from \cref{thm:main}.
\end{theorem}

\begin{proof}
As before, we begin with the bounds in \cref{u/v decay,derbound,rp gain}, which in the inside region translate to
\begin{equation}\label{inbd1}
  \p_{\le m+n} \ls \frac{\ju^{1/2}}{\la t\ra}, \quad \pa\p_{\le m+n} \ls \frac{1}{\la r\ra \la u\ra^{1/2}}, \quad \pmn \ls  \frac{\ju^{1/2 - \xg/2}}{\la t\ra}
\end{equation}
where again for simplicity we use the far left $\pmn$ bound for $\p_1$ and the far right $\pmn$ bound for $\p_3$. 

Let $\xs$ denote the reduced, irrational number mentioned in \cref{smllirr.xs} until stated otherwise. We thus have, using \eqref{inbd1}:
\[
H_1 \lesssim \frac{\ju^{1/2}}{\jr^{2+\sigma}\jt}, \quad \pa_t H_2 \lesssim \frac{1}{\la r\ra^{1+\sigma}\jt\la u\ra^{1/2}}, \quad H_3 \ls \f{\ju^{5(1/2 - \xg/2)}}{\jt^5} \ls \f{\ju^{1/2}}{\jt^{3 + 5\xg/2}}.
\] %

By \eqref{lindcy1} with $\alpha=2+\sigma$, $\beta = 1$, and $\eta=-1/2$, we obtain 
\[
(\p_1)_{\leq m+n} \lesssim \jr^{-1} \ju^{1/2-\sigma} %
\]
	Similarly \eqref{lindcy1} with $\alpha=2+\sigma$, $\beta = 0$, and $\eta=1/2$ yields
\[
(\p_2)_{\leq m+n} \lesssim \jr^{-1} \ju^{1/2-\sigma} %
\]
	Finally, \eqref{lindcy1} with $\alpha=0$, $\beta = 3 + 5\xg/2$, and $\eta=-1/2$ yields 
\[
(\p_3)_{\leq m+n} \lesssim \jr^{-1} \ju^{1/2-5\xg/2}
\]

The three inequalities above give
\[
\p_{\leq m+n} \ls \jr^{-1} \ju^{1/2-\sigma'}, \quad \sigma' := \min(2\xg, \xs).
\]
By \cref{convrsn} we obtain the following improved bounds (by a factor of $\ju^{-\sigma'}$): 
\begin{equation}\label{inbd2}
  \pmn \ls \frac{\ju^{1/2-\sigma'}}{\la t\ra}, \quad \pa\pmn \ls \frac{1}{\la r\ra \la u\ra^{1/2+\sigma'}}.
\end{equation}

We now repeat the iteration, replacing $\eta$ by $\eta+\sigma'$, applying \eqref{lindcy2} and then improving decay in $\{ t/r > 2 \}$ by \cref{convrsn}. The process stops after $\lfloor\frac{1}{2\sigma'}\rfloor$ steps, when \cref{lindcy1,convrsn} yield
\begin{equation}\label{3rdbd'}
  |\p_{\le m+n}| \ls \frac{1}{\la t\ra}, \quad |\pa\p_{\le m+n}| \ls \frac{1}{\la r\ra \la u\ra}.
\end{equation}

At this point we switch to using \eqref{lindcy1der} for $\p_2$, and the iteration process follows the same pattern as in \cref{sec:ext}, with the extra use of \cref{convrsn}.
	Like before in \cref{ext.thm}, we make the final iterate involving the \textit{original} value of $\xs$ from \cref{thm:main}.
\end{proof}

\begin{remark}[Other integral powers $p \in \Z_{\ge2}$]\label{rem:HP.iteration}
For $\p$ solving $P\p = G_{3}$ with
 $$G_{3} = |\p|^{p}\p, \quad p\ge5 \text{ with large initial data}$$
 or
  $$G_{3} = \pminus \p^{p+1}, \quad p\ge5 \text{ with small initial data}$$
 we remark that \textit{if} global existence holds for the large data case (as this is clear for the small data case), then just
by the initial global decay rate \cref{u/v decay} alone, the decay rates in \cref{rp gain} are immediately achieved for some $\xg$ (and hence there is no need to prove \cref{the rp est}), and letting $\p_{3}$ solve $\Box\p_{3} = G_{3}$ as before in \cref{settingup}, the iterations in \cref{sec:ext,sec:int} follow nearly verbatim, with the natural modification that in the end we reach the final decay rate 
$$(\p_{3})_{\le m+n} \ls r\inv \ju^{-(p-1)} , \quad \text{ for } u < -1, \qquad (\p_{3})_{\le m+n} \ls v\inv \ju^{-(p-1)} , \quad \text{ for } u > 1	.$$
Thus 
$$\pm\ls \jv\inv \ju^{ - \min \{ 1 + \sigma,  p - 1  \} }$$

For $\p$ solving $P\p = G_{3}$ with
 $$G_{3} = |\p|^{2}\p$$
with \textit{large initial data}, we remark here that \textit{if} global existence holds and if \cref{rp gain} holds with $\gamma >1$, 
then the iteration also follows \cref{sec:ext,sec:int} nearly verbatim, and we reach the final decay rates
$$(\p_{3})_{\le m+n} \ls r\inv \ju\inv  \text{ for }u < -1	,\qquad (\p_{3})_{\le m+n} \ls v\inv \ju\inv \text{ for }u > 1	$$
and thus
$$\pm\ls \jv\inv \ju^{-1 }.$$
Similarly we have global existence for $P\p = \pminus \p^{3}$ with small initial data, and if we had \cref{rp gain} with $\xg>1$ then this remark holds as well.

For $\phi$ solving $P\phi = G_{3}$ with 
$$G_{3} = \pminus \phi^{4}$$
with small initial data, if we had \cref{rp gain} with $\xg>1/2$ then after setting $\Box\p_{3}=G_{3}$ in \cref{settingup}, the iterations in \cref{sec:ext,sec:int} also hold nearly verbatim and we reach
$$(\p_{3})_{\le m+n} \ls \jv\inv \ju^{-2}.$$
thus 
$$\pm\ls \jv\inv\ju^{-\min\{1+\xs, 2\}}.$$

\end{remark}

\end{document}